%% file: PentaI.tex
\documentclass[11pt]{article}
\usepackage{times} 
\usepackage{a4} 
\pdfoutput=1
\usepackage{latexsym} 
\usepackage{amssymb}
\usepackage{amsmath,amsthm}
\usepackage{url}
\usepackage[pdftex]{graphicx}
\DeclareGraphicsExtensions{.jpg,.pdf,.mps,.png}

\newtheorem{lemma}{Lemma}[section]
\newtheorem{proposition}{Proposition}[section]
\newtheorem{claim}{Claim}
\numberwithin{equation}{section} 
\newcommand{\mymax}[1]{\stackrel{\wedge}{#1}} 
\newcommand{\mymin}[1]{\stackrel{\vee}{#1}} 
\newcommand{\meqdef}{\overset{\mathrm{def}}{=}} 
\newcommand{\bminus}{\boldsymbol{-}} 

\begin{document}
\addtocounter{footnote}{1} \footnotetext{Department of Mathematics and 
Computer Science, University of Catania, Italy.\\ 
E-mail:\mbox{\ }scollo@dmi.unict.it}  
\addtocounter{footnote}{-1} 

\title{An integration of Euler's pentagonal partition} 
\author{Giuseppe Scollo}
\date{draft, version 1\\~\\ 19 September 2010} 

\maketitle

\begin{abstract} 
\noindent 
A recurrent formula is presented, for the enumeration of the compositions of 
positive integers as sums over multisets of positive integers, that closely 
resembles Euler's recurrence based on the pentagonal numbers, but where the 
coefficients result from a discrete integration of Euler's coefficients. 
Both a bijective proof and one based on generating functions show the 
equivalence of the subject recurrences.
\end{abstract}

\thispagestyle{empty}

\input{PentaIs1}

\input{PentaIs2}

\input{PentaIs3}

\input{PentaIs4}

\input{PentaIs5}

\input{PentaIs6}

\input{PentaIs7}

\input{PentaIak}

\medskip
\bibliographystyle{unsrt} 
\bibliography{PentaI} 

\end{document}

%% file: PentaIs1.tex
\section{Introduction} 
\label{s1} 
Euler's pentagonal recurrence for integer partitioning \cite{E1753} may be 
presented as follows. By ancient Greek tradition, pentagonal numbers are 
those of the form $(3m^2-m)/2$. To get all which are needed for Euler's 
recurrence, the range of $m$ is extended to all integers, including the 
negative ones. Euler's coefficients for pentagonal recurrent partitioning, 
then, form the following sequence, indexed by the natural numbers: 
\begin{equation}
e_{_{n}} = (-1)^{k+1} \mbox{\rm\ \ \ \ \ if\ } n=(3k^2\pm k)/2, 
\mbox{\rm\ \ \ \ \ } 
e_{_{n}} = 0 \mbox{\rm\ \ \ \ \ if\ } n \mbox{\rm\ is not pentagonal.} 
\label{eq:1.1} 
\end{equation}
Euler's pentagonal partitioning may then be obtained by the following 
recurrence, having set forth that $p(0) = 1$ and p($n$) = 0 for negative $n$: 
\begin{equation}
p(n) = \mathop{\sum}_{{k>0}} e_{_{k}}p(n-k) 
\label{eq:1.2}
\end{equation}
Now, consider the following sequence of coefficients, which result from a 
discrete integration of Euler's sequence: 
\begin{equation}
f_{_{n}} = \mathop{\sum}_{{0\leq k\leq n}} e_{_{k}} \mbox{\rm\ for\ } n\geq 0. 
\label{eq:1.3}
\end{equation}
It's easy to see that $f_{_{n}} = e_{_{n}}$ iff 
$n=0$\ or $(3m^2-m)/2<n\leq(3m^2+m)/2$\ for some positive $m$, and that, 
just like Euler's coefficients, also those defined by equation (\ref{eq:1.3}) 
are bound to take values in $\{0,\pm 1\}$.  

We claim the following recurrence holds as well. 
\begin{claim} 
The coefficients defined by Equation~(\ref{eq:1.3}) satisfy the recurrence: 
\begin{equation}
p(n) = 1+\mathop{\sum}_{{k>0}} f_{_{k}}p(n-k)
\label{eq:1.4}
\end{equation}
\end{claim} 
The proof of the validity of our claim is deferred until Section\,\ref{s6}, 
however, whereas the forthcoming sections aim at elucidating its combinatorial 
as well as computational roots. 

Let's fix some notation and terminology, for the purposes of the present note: 
\begin{itemize} 
\item 
${\cal P}_n$ : the set of {\em partitions} of natural number $n$, viz.\ the 
multisets of positive integers whose sum is $n$; the elements of a partition 
are referred to as its {\em parts}; 
\item 
${\cal S}_n$ : the set of {\em strict partitions} of $n$, which are those 
where all parts are distinct, {\em i.e.} every part has multiplicity 1;  
\item 
${\cal P}_{n_{_{\!<k}}}, {\cal P}_{n_{_{\!>k}}}$ : the subset of ${\cal P}_n$ 
where every part is constrained to be lower, resp.\ higher than $k$; 
\item 
${\cal P}_{n_{\mymax{k}}}, {\cal P}_{n_{\mymin{k}}}$ : the subset of 
${\cal P}_n$ where $k$ is the largest, resp.\ smallest part; 
\item 
${\cal S}_{n_{_{\!<k}}}, {\cal S}_{n_{_{\!>k}}}$, 
${\cal S}_{n_{\mymax{k}}}, {\cal S}_{n_{\mymin{k}}}$ : 
the similarly defined subsets of ${\cal S}_n$; 
\item 
$p(n), s(n),$\ 
$p_{\!_{<k}}(n), p_{\!_{>k}}(n), p_{\mymax{k}}(n), p_{\mymin{k}}(n),$\ 
$s_{\!_{<k}}(n), s_{\!_{>k}}(n), s_{\mymax{k}}(n), s_{\mymin{k}}(n)$ : 
the cardinality of ${\cal P}_n, {\cal S}_n,$\ 
${\cal P}_{n_{_{\!<k}}}, {\cal P}_{n_{_{\!>k}}},$\ 
${\cal P}_{n_{\mymax{k}}}, {\cal P}_{n_{\mymin{k}}},$\ 
${\cal S}_{n_{_{\!<k}}}, {\cal S}_{n_{_{\!>k}}},$\ 
${\cal S}_{n_{\mymax{k}}}, {\cal S}_{n_{\mymin{k}}}$, 
respectively. 
\end{itemize}

%% file: PentaIs2.tex
\section{Recurrences for integer partitioning} 
\label{s2} 
Several recurrences are known to compute $p(n)$, see {\em e.g.}\ 
\cite{A1983,A1998,P2006,W2000}. Some are {\em direct} recurrences, 
in the sense that only the subject function occurs as a recurrent 
in the recurrence body, {\em e.g.} as it happens with Equation (\ref{eq:1.2}). 
Their implementation by dynamic programming only takes $O(n)$ space to store 
the only-once computed values of the recurrents, for a given input $n$. 
Another well-known direct recurrence for integer partitioning is the 
following, also originating from Euler's investigations \cite{E1783} 
(see \cite{B2006} for a history of Euler's work on the pentagonal number 
theorem): 
\begin{equation} 
p(n) = \frac{1}{n}\mathop{\sum}_{{k\geq 1}}\sigma(k)p(n-k) 
\label{eq:2.1} 
\end{equation}
where $\sigma(k)$ is the sum of the divisors of $k$. This recurrence may be 
obtained by a straightforward manipulation of Euler's generating function 
for $p(n)$: 
\begin{equation} 
\mathop{\prod}_{{j\geq 1}}\frac{1}{1-x^j} = \mathop{\sum}_{{n\geq 0}}p(n)x^n.  
\label{eq:2.2}
\end{equation}
The method that enables one to get a recurrence out of a generating function, 
such as (\ref{eq:2.1}) from (\ref{eq:2.2}), is well-known (see {\em e.g.} 
\cite{W2000}, pp.\ 8--9, for a clear exposition), and we do not deal with it 
now, but we anticipate that we make use of it in Section~\ref{s6.2}, where 
it turns out to be a helpful tool to prove the main claim of this note by 
means of generating functions. By the way, also the claimed Equation 
(\ref{eq:1.4}) is a direct recurrence. 

Besides direct recurrences, several recurrences of a different kind are known 
for integer partitioning; their common character is, of course, that 
they involve the use of an {\it auxiliary} recurrence, that depends on 
additional parameters, most commonly one, such as a recurrence for any of 
the bound-indexed partition functions listed at the end of Section~\ref{s1}. 
Actually, recurrences that depend on additional parameters also find 
applications on their own, {\em e.g.} for computational purposes such as 
those reported by \cite{W2000}, p.\ 13. For the purposes of the present note, 
however, our primary interest is in their use as auxiliary devices, to get 
a closer insight into the combinatorial justification of (usually direct) 
recurrences obtained by other means, such as the analytical manipulation of 
generating functions. A classical, highly relevant example in this respect is 
Franklin's combinatorial proof \cite{F1881} of Euler's pentagonal number 
theorem (see {\em e.g.}\ \cite{Z2003} for a tutorial exposition of Franklin's 
proof). 
More recently, a note by Kevin Brown in his math pages \cite{Brown} 
illustrates an enlightening bijective proof of Euler's pentagonal 
recurrence that, unlike Franklin's proof, doesn't even make use of the 
fact that the generating function of the pentagonal coefficients and 
that of the partition function are reciprocal. 

Because of their composite functional structure, involving mutual recurrence 
between distinct recursive functions, we call {\em composite} recurrences for 
integer partitioning those where auxiliary recurrences occur. Here are a few, 
well-known examples, which turn out to be relevant to the developments in 
the forthcoming sections, together with their combinatorial justification. 
The first example is the composite recurrence adopted in the aforementioned 
note, that makes use of an auxiliary recurrence on $p_{\mymin{k}}(n)$, 
the number of partitions of $n$ with smallest part $k$. This satisfies the 
following equations: 

\begin{equation}
p(n) = p_{\mymin{1}}(n+1) 
\label{eq:2.3}
\end{equation}
\begin{equation}
p_{\mymin{k}}(n) = \mathop{\sum}_{{i\geq k}}p_{\mymin{i}}(n-k) 
\mbox{\rm\ \ if\ } k<n 
\label{eq:2.4}
\end{equation}
\begin{equation}
p_{\mymin{n}}(n) = 1 
\label{eq:2.5}
\end{equation}
\begin{equation}
p_{\mymin{k}}(n) = 0 \mbox{\rm\ \ for\ } k>n 
\label{eq:2.6}
\end{equation}
\begin{equation}
p_{\mymin{k+1}}(n) = p_{\mymin{k}}(n-1) - p_{\mymin{k}}(n-1-k) . 
\label{eq:2.7}
\end{equation}

The evidence of the first equation is immediate; we just point out its r\^ole 
in the reduction of the computation of $p(n)$ to that of the auxiliary 
recurrent partitions, thanks to Equation~(\ref{eq:2.4}). The latter is easily 
justified by considering the effect of the removal of a minimal part from 
each of the partitions in ${\cal P}_{n_{\mymin{k}}}$; one clearly gets the 
set ${\cal P}_{(n\bminus k)_{_{\!>k\bminus 1}}}$, 
whose cardinality may be computed 
by summing up $p_{\mymin{i}}(n-k)$ for all $i\geq k$. These contributions 
may be computed by using Equation~(\ref{eq:2.7}), together with the obvious 
basis provided by Equations~(\ref{eq:2.5}--6). A combinatorial argument for 
Equation~(\ref{eq:2.7}) is obtained by considering the transfer of the negative 
term to the left hand side. Then ${\cal P}_{n-1_{\mymin{k}}}$ may be split 
into two disjoint subsets, viz.\ the partitions where the minimal part has 
multiplicity greater than 1, and those where there's only one minimal part. 
The former are clearly counted by $p_{\mymin{k}}(n-1-k)$, again by considering 
the effect of the removal of a minimal part; the latter are counted by 
$p_{\mymin{k+1}}(n)$, by considering the effect of adding 1 to the (only one) 
minimal part. 

Equation~(\ref{eq:2.7}), deployed as a left-to-right computation rule, warrants 
reduction of the computation of any auxiliary term $p_{\mymin{k+1}}(n)$ to 
terms with minimal part 1 and lower $n$, hence to contributions to a direct 
recurrence for $p(n)$. This turns out to be Euler's pentagonal 
recurrence~(\ref{eq:1.2}), details may be found in the aforementioned note. 
A relevant feature of the computational reduction displayed above, is 
the {\em difference} form of the right hand side of Equation~(\ref{eq:2.7}). 
This tells why may it happen that most of the contributions yield a null 
result, which must be the case to get a recurrence with so many null 
coefficients as Euler's one. Such a feature is not enjoyed by other composite 
recurrences, such as the following one, making use of an auxiliary recurrence 
on $p_{\!_{<k}}(n)$, which counts partitions with upper-bounded parts. 
The basic idea is to recursively split ${\cal P}_n$ into two disjoint subsets 
of partitions: those with maximal part $k$, and those where every part is 
lower than $k$, with $k$ ranging from $n$ down to 2. The top-level split 
enables one to get $p(n)$ by recursively computing $p_{\!_{<k}}(n)$, using 
the following equations: 

\begin{equation}
p(n) = 1 \mbox{\rm\ \ for\ } 0\leq n\leq 1 
\label{eq:2.8} 
\end{equation}
\begin{equation}
p(n) = 1 + p_{\!_{<n}}(n) \mbox{\rm\ \ for all\ } n\geq 2 
\label{eq:2.9} 
\end{equation}
\begin{equation}
p_{\!_{<2}}(n) = 1 \mbox{\rm\ \ for all\ } n\geq 0 
\label{eq:2.10}
\end{equation}
\begin{equation}
p_{\!_{<k+1}}(n) = \mathop{\sum}_{m=0}^{\lfloor n/k\rfloor}p_{\!_{<k}}(n-mk) %
\mbox{\rm\ \ \ for\ \ } 2\leq k< n  
\label{eq:2.11}
\end{equation}
\begin{equation}
p_{\!_{<k}}(n) = p(n) \mbox{\rm\ \ if\ \ } k>n 
\label{eq:2.12}
\end{equation}

The combinatorial evidence of these equations needs little explanation; 
it may be useful to point out that Equation~(\ref{eq:2.11}) splits 
${\cal P}_{n_{_{\!<k+1}}}$ into $\lfloor n/k\rfloor+1$ pairwise disjoint 
subsets, according to the multiplicity $m$ of the maximum allowed part 
(viz.\ $k$) as maximal part in the partition, for the given upper bound on it 
(for $m=0$\ one thus gets the partitions where maximal parts are stricly lower 
than $k$). 
On the other hand, neither is this recurrence computationally convenient (its 
dynamic programming implementation takes $O(n^2)$ space to store the computed 
recurrents), nor does it immediately lend itself to reduction to a direct 
recurrence where a significant subset of the coefficients would be null, since 
the deployment of its equations as computational rules, unlike the previous 
case, features no difference of recurrents in the right-hand-side. However, 
it does offer a good basis for further combinatorial reasoning, that leads to 
a different, composite recurrence which enjoys this property, as it is shown 
in Section~\ref{s5}.

%% file: PentaIs3.tex
\section{Auxiliary reductions in composite recurrent partitioning} 
\label{s3} 

A common feature of composite recurrences of interest in this note, is 
their ability to inductively reduce the computation of terms of form 
$p_{\!_{\diamond k}}(n)$, where $\diamond$ is a generic designator of the 
type of constraint that is imposed over partitions counted by the auxiliary 
recurrence, to terms of similar form, but with such values of the $k,n$ 
pair that they are equated, by the {\em relative} inductive basis of the 
given recurrence, to terms $p(n-j)$ of the underlying direct recurrence, 
thus for $1\leq j\leq n$. The {\em relative} qualification is precisely meant 
to say here that, by induction on the auxiliary parameter, all auxiliary terms 
reduce to terms of the other family, so, the relative basis of the induction 
does not consist of a computational assignment of values to the auxiliary 
terms which have minimal values of the induction parameter, it rather consists 
of their immediate reduction to nonauxiliary terms. 

By the way, in Section~\ref{s1} we only introduced the four types of 
constraints which are relevant to the present note, but many other types may 
well deserve interest in other contexts. For example, it is easy to devise 
constructions of direct recurrences by induction on the cardinality of  
partitions, and hereby, if so wished, to consider minima and maxima thereof, 
upper and lower bounds thereupon, etc.. Constraints may be combined even 
further, thereby giving rise to auxiliary recurrences with more than one 
parameter. Alternative auxiliary recurrences do not always deliver different 
outcomes, though. For example, it is well known that induction on the 
cardinality of partitions, or on upperbounds thereupon, is equivalent to 
induction on maximal parts, or on upperbounds thereupon, respectively. 
This is immediately seen by transposing the Ferrers diagrams which represent 
partitions \cite{W2000}. This very fact, however, also shows that such 
equivalences do not hold for strict partitions, since transposition of 
Ferrers diagrams does not preserve strictness. 

If the construction of a direct recurrence is meant to be the purpose 
of the composite one, then a greater interest arises in the 
 {\em auxiliary reductions} produced by the composite recurrence, 
in order to evaluate their contributions to coefficients of the target 
direct recurrence. More precisely, this is formalized as follows. 
For the sake of simplicity, only composite recurrences with one auxiliary 
parameter are considered. Generalization to the multiparameter case is 
straightforward, but not needed for the purposes of this note. 

The following concept proves useful to the forthcoming formalization. Recall 
that a {\em rewrite rule} is a pair of terms with variables, usually written 
in the form $t_1\,\rightarrow\,t_2$, such that 1) $t_1$ is not a variable, and 
2) every variable which occurs in $t_2$ also  occurs in $t_1$. 
A rewrite rule may be extended with a {\em domain condition}, viz.\ a predicate 
with variables which occur in $t_1$, that specifies the rule applicability 
domain. Rewrite rules with domain conditions may be written in the form 
$[d] \; t_1\,\rightarrow\,t_2$, with $d$ the domain condition. 

Rule instantiation, being a syntactic operation, is not constrained by domain 
conditions---but a rule instance results from applying a substitution to all 
variable occurrences in all rule constituents, domain condition included. 
Domain conditions rather affect the definition of {\em ground rewriting system} 
generated by a set ${\cal R}$ of rewrite rules with domain conditions. This is 
the set of ground rewrite rules, viz.\ rewrite rules with neither variables 
nor domain conditions, that are $\sigma$-instances of some rewrite rule $r$ 
by a closed substitution $\sigma$ such that 1) there is a rule $[d] \; r$ in 
${\cal R}$, and 2) the interpretation of ground predicate $\sigma d$ holds. 
Let ${\cal R}_{\mathbf{g}}$ denote the ground rewriting system generated by 
the set ${\cal R}$ of rewrite rules with domain conditions. 

For a given target partition-counting function and auxiliary recurrence 
function, the set of {\em primary recurrence atoms} is defined to consist of 
the terms of form $P(u)$, where $P$ is a generic designator of the target 
function (which may be $p$, $s$, or any other which may be of interest), 
while the set of {\em auxiliary recurrence atoms} consists of the terms of 
form $A(u,v)$, where $A$ is a generic designator of the auxiliary recurrence 
function, with $u, v$\ arithmetic terms, possibly with variables, that 
may only take integer values. 
A recurrence {\em literal} is either a recurrence atom or the product of a 
recurrence atom by an integer arithmetic term; the literal is either primary 
or auxiliary depending on the similar qualification of its constituent atom. 

Assume now we are given a finite set ${\cal R}$ of rewrite rules with domain 
conditions, on these literals extended with {\em additive} arithmetic terms 
built upon literals and integer arithmetic terms. 
Without too much loss of generality, assume that for every rule 
$[d] \; t_1\,\rightarrow\,t_2$ in ${\cal R}$, $t_1$ is a recurrence atom with 
only variables as proper subterms, while $t_2$ is an additive arithmetic term 
on recurrence literals and integer arithmetic terms; this term is assumed to 
be in additive arithmetic normal form, here defined as an indexed sum of 
recurrence literals, where each atom may occur in at most one literal, plus 
at most one standalone arithmetic term. Index bounds may well be integer 
arithmetic terms. Indexing may be implicit, whenever binary additive operators, 
rather than an explicitly indexed sum operator, constitute the sum; in this 
case the index assigned to each literal is identified with the sequential 
position of the recurrence literal in the sum term. 

Further, assume that ${\cal R}$ contains one or more rules where the left 
hand side term is a primary atom while one or more auxiliary atoms and no 
primary atom occur in the right hand side term, as well as one or more rules 
where the left hand side term is an auxiliary atom while no auxiliary atom 
occurs in the right hand side term. The former are referred to as 
 {\em startup rules}, the latter as {\em termination rules}, while all other 
rules are assumed to belong to either of the following categories: 
 {\em primary rules}, where the left hand side term is a primary atom while 
no auxiliary atom occurs in the right hand side term, and 
 {\em auxiliary rules}, where the left hand side term is an auxiliary atom 
while one or more auxiliary atoms and no primary atom occur in the right hand 
side term. We thus rule out only rules (pun intended) where both primary atoms 
and auxiliary atoms occur in the right hand side term. 

Rules may be (uniquely) named, and thus be put in the following general form, 
for each of the four rule types just introduced, where $r$ is the rule name 
(which may also be made use of to designate the rule itself). Please note that 
the summation index upperbound is generally allowed to be lower than the 
lowerbound, in which case the summation is null. However, in those cases where 
the upperbound is not allowed to be lower than the lowerbound, this constraint 
is specified to the right of the index upperbound in the summation. 
The lowerbound 1 is taken in the following forms, with no loss of generality 
since one may always meet this assumption by an 
 index substitution in a given summation. 
%
\begin{subequations} 
\begin{align} 
\mathrm{primary\ } r &: \quad [d_r(n)] \quad \!\! \qquad \qquad 
P(n) \quad \rightarrow \quad 
t_{r_0} + \mathop{\sum}_{1\leq i\leq n_r} t_{r_i} P(u_{r_i}), 
\label{eq:3.1a}\\ 
\mathrm{startup\ } r &: \quad [d_r(n)] \quad \!\! \qquad \qquad 
P(n) \quad \rightarrow \quad 
t_{r_0} + \mathop{\sum}_{1\leq i\leq n_r\geq 1} t_{r_i} A(u_{r_i},v_{r_i}), 
\label{eq:3.1b}\\ 
\mathrm{auxiliary\ } r &: \quad [d_r(n,k)] \quad \!\! \qquad 
A(n,k) \quad \rightarrow \quad 
t_{r_0} + \mathop{\sum}_{1\leq i\leq n_r\geq 1} t_{r_i} A(u_{r_i},v_{r_i}), 
\label{eq:3.1c}\\ 
\mathrm{termination\ } r &: \quad [d_r(n,k)] \quad \!\! \qquad 
A(n,k) \quad \rightarrow \quad 
t_{r_0} + \mathop{\sum}_{1\leq i\leq n_r} t_{r_i} P(u_{r_i}), 
\label{eq:3.1d} 
\end{align} 
\end{subequations} 
where $n,k$ are variables ranging over the integers, 
$t_{r_i},u_{r_i},v_{r_i},n_r$ denote integer arithmetic terms on these 
variables as well as the bound variable $i$, for $1\leq i\leq n_r$, but not on 
$k$ in primary and startup rules unless it is the bound variable of the 
summation in the scope of which those terms occur, $t_{r_0}$ optional, but 
mandatory if $n_r = 0$.  

Let ${\cal R}_P, {\cal R}_I, {\cal R}_A$, and ${\cal R}_T$, denote the subsets 
of ${\cal R}$ that consist of the primary, startup, auxiliary, and termination 
rules in ${\cal R}$, respectively. Furthermore, let 
${\cal R}_1 = {\cal R}_P \cup {\cal R}_I$, and 
${\cal R}_2 = {\cal R}_A \cup {\cal R}_T$.  

Let $\tilde{t}$\ denote the value of arithmetic ground term $t$. 
A rewrite system ${\cal R}$ composed of rules of the form displayed above, 
is {\em unitary} if every ground rule $r \in {\cal R}_{\mathbf{g}}$ it 
generates, satisfies the following condition: 
\begin{equation} 
1 \leq i \leq n_r \Rightarrow 
((-1 \leq \tilde{t}_{r_i} \leq 1\}) 
 \wedge 
 ((j > 0 
   \wedge
   (\tilde{u}_{r_i}=\tilde{u}_{r_j}) 
   \wedge 
   (\tilde{v}_{r_i}=\tilde{v}_{r_j}) 
  )\Rightarrow i=j 
)) . 
\label{eq:3.2} 
\end{equation} 
We henceforth assume to deal with unitary rewrite systems, since they suffice 
to the purposes of the present note, although much of the forthcoming work may 
be extended to nonunitary rewrite systems by a straightforward generalization. 

Furthermore, an {\em orthogonality} requirement is put on ${\cal R}$, that 
bears some resemblance with the analogous, syntactic property as defined 
for term rewriting systems, but in the present context it generally depends 
on the interpretation of (ground) domain conditions. Briefly, the aim is 
to make sure that every ground rule in ${\cal R}_{\mathbf{g}}$ may be traced 
back to only one rule in ${\cal R}$. This is warranted by the requirement 
that domain conditions of rules in ${\cal R}$ specify pairwise disjoint sets 
of ground instances whose interpretation holds, that is, for every closed  
substitution $\sigma$ and every pair $r, r^{\prime}$ of rules in ${\cal R}_i$, 
with $1\leq i\leq 2$, letting $d_r, d_{r^{\prime}}$ denote their respective 
domain conditions (with concise, but somewhat cavalier notation, hopefully 
forgiven by the learned reader): 
\begin{equation}
(\sigma d_r \wedge \sigma d_{r^{\prime}}) \Rightarrow\ r = r^{\prime} 
\label{eq:3.3} 
\end{equation} 

The additive shape of terms which are assumed to form the right hand side 
of auxiliary rules, together with the orthogonality and unitarity assumptions 
enable the following geometric interpretation of their ground instances, 
hereafter termed {\em (auxiliary) ground rules}, for brevity. 
The {\em parallel reduction graph} of auxiliary ground rules may be construed 
as a directed acyclic graph (DAG) with labelled edges and (also labelled) 
vertices in the discrete Cartesian plane. Ground auxiliary atoms $A(u,v)$ are 
interpreted as points $(\tilde{u},\tilde{v})$ in the plane. 
Since ${\cal R}$ is unitary, each auxiliary ground rule may be put in the 
following form 
\begin{equation} 
r : A(u_{r_0},v_{r_0})\rightarrow 
    t_r + \mathop{\sum}_{1\leq i\leq n_r\geq 1} s_{r_i} A(u_{r_i},v_{r_i}), 
\label{eq:3.5}
\end{equation}
with $n_r, u_{r_i}, v_{r_i}$ arithmetic ground terms, for 
$0\leq i\leq \tilde{n}_r$, and 
$s_{r_i}\in\{\pm 1\}$, for $1\leq i\leq \tilde{n}_r$, 
and contributes a fan of edges to the DAG construction, all outgoing from 
$(\tilde{u}_{r_0},\tilde{v}_{r_0})$, and each respectively incoming to 
$(\tilde{u}_{r_i},\tilde{v}_{r_i})$; moreover, ground rule $r$ contributes 
a {\em sign label} $s_{r_i}$\ to the edge incoming to target vertex 
$(\tilde{u}_{r_i},\tilde{v}_{r_i})$ in the fan, for $1\leq i\leq\tilde{n}_r$, 
and, whenever a standalone arithmetic term $t_r$ occurs in the right hand side 
of the rule, this contributes an integer {\em constant label} $\tilde{t}_r$ to 
the source vertex of the fan, $(\tilde{u}_{r_0},\tilde{v}_{r_0})$, otherwise 
labelled by the default label 0 (omitted when drawing the DAG). 
Thanks to the orthogonality requirement, every vertex and every edge of the 
DAG have label provided by a unique rule, while the unitarity assumption 
warrants two-valuedness of edge labels. This assumption may be relaxed by 
taking integers as edge labels. 

Before addressing the interpretation of startup and termination rules in the 
discrete Cartesian geometry, it is convenient to check whether the assumptions 
made so far may be met in cases of interest. To this purpose, first, consider 
the composite recurrence defined by equations~(\ref{eq:2.3})--(\ref{eq:2.7}). 
As they stand, their left-to-right reading as rewrite rules does not comply 
with the orthogonality requirement on the auxiliary rewrite system; {\em e.g.}, 
rules corresponding to equations~(\ref{eq:2.4}) and (\ref{eq:2.7}) overlap. 
One may get an equivalent set of equations (with domain conditions), whose 
left-to-right reading complies with that requirement as well as all other 
assumptions made so far, as follows. Take the equation obtained by transitivity 
from equation~(\ref{eq:2.3}) and the [$k\mapsto 1, n\mapsto n+1$]-instance 
of equation~(\ref{eq:2.4}), and put upperbound $n$ and lowerbound 1 to the 
summation index, thanks to equation \ref{eq:2.6}. A basis equation for 
$p(0)$ is separately needed. Equation~(\ref{eq:2.7}) should to be limited 
to its ($1\leq k< n-1$)-instances, not to overlap with 
equations~(\ref{eq:2.5}--6), but we may as well raise by 1 the upperbound 
on $k$ and dispose of equation ~(\ref{eq:2.5}). 
Equation~(\ref{eq:2.3}) is turned into a rewrite rule by right-to-left reading, 
Equation~(\ref{eq:2.4}) may be safely disposed of. 
All this results in the following rewriting system with domain conditions: 
\begin{subequations}
\begin{align}
[n = 0] \; \; \qquad p(n) \rightarrow & \; 1 
\label{eq:3.6a}\\ 
[n > 0] \; \; \qquad p(n) \rightarrow & 
\mathop{\sum}_{{1\leq i\leq n}}p_{\mymin{i}}(n) 
\label{eq:3.6b}\\ 
[n > 0 \wedge k = 1] \qquad 
p_{\mymin{k}}(n) \rightarrow & \; p(n-1) 
\label{eq:3.6c}\\ 
[2\leq k \leq n] \qquad p_{\mymin{k}}(n) \rightarrow & \; 
p_{\mymin{k-1}}(n-1) - p_{\mymin{k-1}}(n-k) 
\label{eq:3.6d}\\ 
[k > n] \qquad p_{\mymin{k}}(n) \rightarrow & \; 0 
\label{eq:3.6e} 
\end{align}
\end{subequations}
Clearly, this system is unitary and orthogonal. It is composed of one primary 
rule~(\ref{eq:3.6a}), one startup rule~(\ref{eq:3.6b}), two termination 
rules~(\ref{eq:3.6c}), (\ref{eq:3.6e}), and one auxiliary rule~(\ref{eq:3.6d}). 

\begin{figure}[hbt]
\centering
\includegraphics[width=14.9cm]{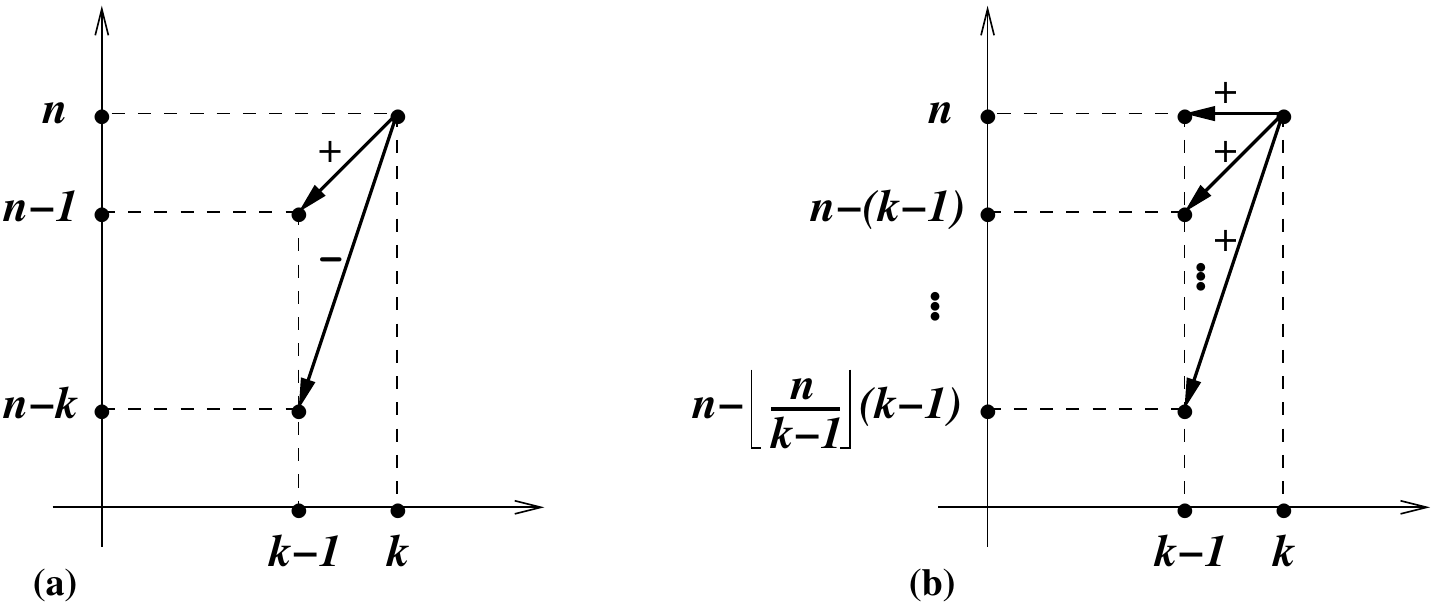}
\caption{Auxiliary ground rules for unitary recurrences}
\label{fig:1}
\end{figure}
 
Figure~\ref{fig:1}(a) displays the fan which represents a generic ground 
instance of rule~(\ref{eq:3.6d}) in the parallel reduction DAG, where the first 
coordinate is on the vertical axis (this perhaps unusual choice is motivated 
in Section~\ref{s4}).  

As a second check, consider the composite recurrence defined by 
equations~(\ref{eq:2.8}--12). Unlike in the first case, their left-to-right 
reading as rewrite rules fully complies with all of the assumptions made 
so far about the rewriting systems which are of interest here, but for a 
straightforward introduction of fairly obvious domain conditions in order to 
turn the left hand side atom of equations~(\ref{eq:2.10}--11) into an atom 
whose subscript subterm is a variable. 
We need not reproduce the rewrite rules explicitly, but we just understand 
that they are those equations with left-to-right orientation as rewrite rules, 
and with the easily defined domain conditions as mentioned above. 
This rewriting system is composed of one primary rule~(\ref{eq:2.8}), one 
startup rule~(\ref{eq:2.9}), two termination rules~(\ref{eq:2.10}), 
(\ref{eq:2.12}), and one auxiliary rule~(\ref{eq:2.11}). 
Figure~\ref{fig:1}(b) displays the fan which represents the interpretation of 
a generic ground instance of this rule in the parallel reduction DAG.

%% file: PentaIs4.tex
\section{Getting direct recurrences out of composite ones} 
\label{s4} 

What r\^ole do primary atoms and nonauxiliary rules play with respect to the 
parallel reduction DAG introduced in the previous section? The answer to this 
question will prove straightforward once the target of the representation 
under development is formalized. 

Recall, the purpose is as per title of this section, thus it entails that 
primary rules fit the purpose as they stand, therefore one need not do anything 
with them, but to include them as equations (with domain conditions, this is 
henceforth understood) in the set of equations forming the aimed at direct 
recurrence. An additional bit of information which may be extracted from the 
primary rules is the characterization of the subset of the partitioning domain 
(the nonnegative integers, in the subject case) which they apply to, so that 
its complement in the subject domain may be taken as the domain of the as yet 
to be discovered part of the target recurrence. But this bit is quite a 
redundant one, since the orthogonality requirement implies the domain of 
startup rules is included therein, and it actually coincides with it, if the 
composite recurrence is to define a total function over the subject domain, 
viz.\ the composite recurrence is {\em complete}, which property is henceforth 
assumed. 

In view of the forthcoming formalization, let $D(n)$ be a predicate over 
the integers that characterizes the domain of the target direct recurrence, 
excluding the subdomain covered by the primary rules. Then, by orthogonality 
of ${\cal R}$ and completeness of the composite recurrence, the family 
of subsets that are characterized by predicates in the 
$(d_r(n) \; | \; r\in{\cal R}_I)$ family, partitions the set characterized 
by $D(n)$. One may conceive to design the as yet unknown part of the 
target direct recurrence as a set of ${\cal R}_I$-indexed equations, one 
for each $r\in {\cal R}_I$, that thus bijectively correspond to the rules 
in ${\cal R}_I$. This justifies the choice of naming each of the subject 
equations with the same name as the corresponding startup rule, with no 
danger of confusion, thanks to the different syntactic shapes of equations 
and rewrite rules. Now, if the inductive nature of recurrences as function 
definitions is taken into account, then it is easy to realize the convenience 
of giving the following form to the part of the aimed at direct recurrence 
that does not come from primary rules. This part will consist of one $r$-named 
equation for each startup rule in ${\cal R}_I$, the $r$-named one, that is 
designed to be of the following form: 
\begin{equation} 
r : \qquad d_r(n) \quad \rightarrow \quad 
    P(n) = c_{r_0} + \mathop{\sum}_{1\leq j\leq n} c_{r_j} P(n-j), 
\label{eq:4.1}
\end{equation} 
where the {\em coefficient terms} 
$(c_{r_j} \; | \; 0\leq j\leq n,\; r\in{\cal R}_I)$ 
are the whole and essential subject of the design under consideration. Each 
startup rule is thus meant to eventually result in a map $c_r$, that is of 
type $\mathbb{N}\!\rightarrow\!\mathbb{Z}$, with $\mathbb{N}$ the nonnegative 
integers, in the fairly frequent case that the coefficients terms 
in equation~\ref{eq:4.1} are {\em constants}, viz.\ only depend on $j$, 
not on $n$, otherwise it is of type 
$\mathbb{N}\!\rightarrow\!\mathbb{N}\!\rightarrow\!\mathbb{Z}$. 

Now, startup rules have found a way to the target, but not yet one to the 
method, that is to say, to the parallel reductions DAG. To this end, consider 
the lower dimensionality of primary atoms with respect to that of auxiliary 
atoms. Ground instances of auxiliary atoms are interpreted as points of the 
discrete Cartesian space where the DAG lives in, thus it is fairly obvious that 
ground instances of primary atoms be interpreted as points of an isomorphic 
image of a unidimensional subspace of the subject space, viz.\ points of the 
coordinate axis which hosts the interpretation image of their corresponding 
projection in auxiliary atoms. Since the present target is to construct 
a family of integer maps $(c_r \; | \; r \in {\cal R}_I)$, it is fairly 
natural to take a second copy of the discrete Cartesian plane to host the 
representation of the target maps, and to interpret left hand side atoms 
of ground startup rules as points on the {\em second} coordinate axis of this 
plane, which is henceforth referred to as the {\em primary plane}, where the 
coefficient functions of the target primary recurrence are sought for. 
The former plane may be qualified as the {\em auxiliary plane}, and this is 
taken as the default plane, thus unless otherwise specified. 
The reason for the choice of the second coordinate axis to represent primary 
atoms relates to the aforementioned target, since it makes the first coordinate 
axis available to represent the values $j$ of the $c_r$ index (for each 
$r \in {\cal R}_I$), with the second coordinate representing the value of 
$c_{r_j}$ in the constant case (otherwise a third dimension is needed, of 
course). 

Similarly to the instantiation of form (\ref{eq:3.1c}) to form (\ref{eq:3.5}) 
for ground auxiliary rules under the unitarity assumption, ground startup 
rules are of the form 
\begin{equation} 
r : P(u_{r_0})\rightarrow 
    t_r + \mathop{\sum}_{1\leq i\leq n_r\geq 1} s_{r_i} A(u_{r_i},v_{r_i}), 
\label{eq:4.2} 
\end{equation} 
under the same assumption. Consistently with the representation of auxiliary 
ground rules in the DAG, each ground startup rule $r$ contributes a fan of 
edges to the extended DAG construction (now spanning over two planes), all 
outgoing from source vertex $(0, \tilde{u}_{r_0})$ in the primary plane, 
and each respectively incoming to $(\tilde{u}_{r_i},\tilde{v}_{r_i})$. 
Moreover, ground rule $r$ contributes the sign label $s_{r_i}$\ to the edge 
incoming to target vertex$ (\tilde{u}_{r_i},\tilde{v}_{r_i})$ in the fan, 
for $1\leq i\leq n_r$, and, whenever a standalone arithmetic term $t_r$ occurs 
in the right hand side of the rule, this contributes the constant label 
$\tilde{t}_r$ to the source vertex of the fan, $(0, \tilde{u}_{r_0})$, 
otherwise labelled by the default label 0 (omitted in DAG drawing). 

Finally, as it may be expected, each ground termination rule construes a fan 
of edges in the extended DAG, with source vertex in the auxiliary plane 
and target vertices in the first coordinate axis of the primary plane. 
Under the unitarity assumption, ground termination rules are of the form 
\begin{equation} 
r : A(u_{r_0},v_{r_0})\rightarrow 
    t_r + \mathop{\sum}_{1\leq i\leq n_r} s_{r_i} P(u_{r_i}), 
\label{eq:4.3} 
\end{equation} 
and the fan of edges construed by such a rule has source vertex 
$(\tilde{u}_{r_0},\tilde{v}_{r_0})$ with constant label $\tilde{t}_r$, and 
an edge for each ground primary atom $P(u_{r_i})$, 
with target vertex $(\tilde{u}_{r_i},0)$ in the target plane and sign label 
$s_{r_i}$. 

A few conclusions may be easily drawn now, about values of the coefficients 
of the target primary rules (\ref{eq:4.1}), starting with  the ``constant'' 
coefficient $c_{r_0}$. This is only a constant relatively to recurrence, 
for it generally designates a function that returns an integer for each 
value of $n$. It may well happen that this function is a constant one, 
as it often does happen, but that's not necessarily the case in general. 
For each value $\tilde{n}$ of $n$ that satisfies $d_r(n)$, consider 
the subgraph of the extended parallel reduction DAG that is rooted at 
$(0,\tilde{n})$ in the primary plane, and where every terminating path is 
considered up to the point which represents the left hand side atom of the 
ground termination rule which is eventually applied for that path. 
Each vertex in the subgraph contributes the value of its constant label, 
multiplied by the sign product of the edge labels along the path leading 
from the subgraph root to it, to $c_{r_0}(\tilde{n})$. The value of this 
coefficient for the chosen $\tilde{n}$ thus results from the sum of all 
these contributions, including those made by termination rules. 

The sign labels of edges in the fans construed by ground termination rules 
play no r\^ole in computing $c_{r_0}(\tilde{n})$, but they do play one in 
the computation of other coefficients in the target recurrence, for the given 
$\tilde{n}$. In this case, terminating paths in the aforementioned subgraph 
of the DAG are to be considered up to their termination point proper, in the 
primary plane, and each of them contributes a unitary summand, with sign given 
by the sign product of all edge labels along the path, to the value of 
coefficient $c_{r_j}(\tilde{n})$ in the target $r$-equation, where $r$ is 
the startup ground rule which construes the subgraph of the DAG where the 
path lives in, and 
\begin{equation}
j = \tilde{n} - \tilde{u}_{r^{\prime}_i} , 
\label{eq:4.4}
\end{equation}
where $r^{\prime}$\ is the ground termination rule which construes the edge 
fan where the terminal edge of the path is found, and finally $i$ the index 
of its right hand side atom that construes that edge. 

The first conclusions just drawn about the construction of the coefficient 
maps $(c_r \; | \; r \in {\cal R}_I)$, especially equation~(\ref{eq:4.4}),  
invite further analysis. Let $\mathbb{A}_T$\ denote the {\em terminal region} 
of the auxiliary plane, that consists of those points which satisfy the domain 
condition of some termination rule, while let $\mathbb{A}_I$\ denote the 
 {\em startup region}, consisting of those points which interpret an atom 
occurrence in the right hand side of some ground startup rule.  
Since $\tilde{u}_{r^{\prime}_i}$ only depends on 
$\tilde{u}_{r^{\prime}_0}, \tilde{v}_{r^{\prime}_0}$, and the index $i$, the 
$[0\leq i\leq \tilde{n}_{r^{\prime}}]$-indexed family of coefficient 
contributions coming from each point of $\mathbb{A}_T$ 
only depends on the terminal point coordinates, say $(n_t,k_t)$, thanks 
to uniqueness of the applicable rule $r^{\prime}$ at each of its points. 
$\mathbb{A}_T$ may thus be quotiented by an equivalence relation $\theta$ 
that partitions it by identity of coefficient contribution families. 

Now, further progress may be made backwards, by wondering (or wandering, 
in a sort of crabwise strategy, or pondering) where do those contributions 
come from. Each point in $\mathbb{A}_T$ is the termination point of a set 
of paths in the auxiliary plane, each path starting at some point 
$(n_0,k_0) \in \mathbb{A}_I$, and further determined by the sequence of 
indices $(i_m \; | \; 1 \leq m \leq l)$, with $l$ the path length, viz.\ 
the number of edges it consists of, where $i_m$ is the edge index in 
the edge fan construed by the reduction rule which the $m$-th edge in 
the path may be ascribed to. A one-to-one correspondence is thus established 
between terminating paths in the auxiliary plane and pairs 
$(\vec{a_0},\vec{\imath})$, with $\vec{a_0}\!=\!(n_0,k_0)\!\in\!\mathbb{A}_I$ 
and $\vec{\imath} = (i_m \; | \; 1\leq m\leq l)$, with $l$ the path length. 

Finally, as last wandering step, the $\vec{a_0}$ startup point component of 
each of the aforementioned pairs may be replaced by an edge index component 
$i_0$, similar to the constituents of the rest of the structure, that is the 
index of the edge, in the edge fan construed by the startup rule, that connects 
point $\tilde{n}$\ in the second coordinate of the primary plane to the 
auxiliary path startup point $\vec{a_0}$ as target vertex. 

As a matter of fact, a reconsideration of the previous treatment of 
termination rules in the analysis outlined above, suggests the possibility 
to extend the subject paths one step beyond the auxiliary Cartesian plane, 
to include the final edge, ascribed to each of them by a corresponding index 
in the termination rule, that connects the auxiliary path terminal point 
to final point $j$ in the first coordinate axis of the primary plane. 
The full path length is thus $l+2$, where $l$ is that of its part in the 
auxiliary Cartesian plane. It is often the case that, for each $\tilde{n}$\ 
in the second coordinate axis of the primary plane, the set 
$\{\;\tilde{n}_r \;|\; r\in\{r_I(\tilde{n})\}\cup{\cal R_A}\cup{\cal R_T}\;\}$, 
is bounded, where ${r_I(\tilde{n})}$ is the unique startup rule such that 
$\tilde{n}$ satisfies its domain condition. In such a case, letting 
$a_{\mathrm{max}}(\tilde{n})$ denote the maximum element of this set, 
every path may be uniquely coded, for the given $\tilde{n}$, as a word 
of length $l+2$ over an alphabet of cardinality $a_{\mathrm{max}}(\tilde{n})$. 
The coding may be useful if one manages to define two functions on its image: 
a {\em (primary index) valuation} function, returning the coefficient index 
$j$, in the target primary recurrence, which the coded path gives a unitary 
contribution to, and a {\em polarity} function, returning the sign of that 
contribution. 

In practice, it is often sufficient to restrict the analysis outlined so far 
to auxiliary paths only, {\em e.g.} whenever termination rules construct a 
single-edge fan, startup rules label with the same sign all edges in their 
edge fan, and coding of startup points is easily combined with edge sequence 
coding in such a manner that the definition of the two aforementioned functions 
proves straightforward. This is what happens, for instance, in the 
combinatorial proof of Euler's pentagonal recurrence for integer partition 
that is exposed in \cite{Brown}, and it happens as well, albeit with a 
different coding, in the recurrence for the same function that is going to 
be presented next. The valuation equivalence partitioning of the terminal 
region proves useful in this case, and that was the reason to introduce it 
in the first place. 
Nonetheless, it seemed useful to outline a more general framework, that could 
support combinatorial reasoning to solve the target problem, as per title of 
this section, also when the problem instance does not meet conditions, such as 
those listed above, which allow the aforementioned restriction of the analysis.

%% file: PentaIs5.tex
\section{A recurrence for integer partitioning based on maximal parts} 
\label{s5} 

As pointed out at the end of Section~\ref{s2}, the composite recurrence for 
integer partitioning provided by equations~(\ref{eq:2.8}--12) does not lend 
itself to reduction to a direct recurrence with (several) null coefficients, 
because of its lack of difference of recurrents in its right hand side terms. 
However, one may take it as a basis to design another composite recurrence 
which enjoys this feature. We give it the form of a unitary system of rewrite 
rules with domain conditions from the outset. 

The partitioning indexed by maximal parts, $p_{\mymax{k}}(n)$, is taken as 
auxiliary recurrence. The following identity between the two subject 
auxiliary recurrences is plain: 
\begin{equation} 
p_{\!_{<k}}(n) = \mathop{\sum}_{1\leq j < k} p_{\mymax{j}}(n) 
\label{eq:5.1} 
\end{equation} 
The following, immediately evident properties of the latter partitioning prove 
useful: 
\begin{equation} 
p_{\mymax{1}}(n) = p_{\mymax{n}}(n) = p_{\mymax{n-1}}(n) = 1, 
\label{eq:5.2} 
\end{equation} 
as well as the fact that, for $k>n/2$, partitions in ${\cal P}_{n_{\mymax{k}}}$ 
have only one maximal part. Some further combinatorial reasoning justifies the 
following rewrite rules for the subject composite recurrence: 
\begin{subequations} 
\begin{align} 
[n \geq 0] \qquad \,\; p(n) \rightarrow & \; 
1 + \mathop{\sum}_{2\leq k \leq n} p_{\mymax{k}}(n) 
\label{eq:5.3a}\\ 
[\max(2,n/2) \leq k \leq n] \qquad p_{\mymax{k}}(n) \rightarrow & \; p(n-k) 
\label{eq:5.3b}\\ 
[2 \leq k < n/2] \qquad p_{\mymax{k}}(n) \rightarrow & \; 
p_{\mymax{k+1}}(n+1) - p_{\mymax{k+1}}(n-k) 
\label{eq:5.3c} 
\end{align} 
\end{subequations} 
First, the reason for the 1 in equation~(\ref{eq:5.3a}), rather than 0 while 
taking 1 as lowerbound of the summation index, is that by the latter choice 
one would have to put $[n>0]$ as domain condition, and $p(0)\,\rightarrow\,1$ 
should then also be specified as a separate rule. Both choices yield the 
correct rewriting for $p(1)$, under the usual convention that summation is 
null when the index has upperbound lower than lowerbound, but the latter 
choice would produce the erroneous reduction $p(0)\,\rightarrow\,0$, were the 
rule domain not be restricted to the positive integers. Note, however, that 
because of the lower bound on the summation index, rule~(\ref{eq:5.3a}) has 
two {\em primary ground instances}, for $0 \leq n \leq 1$, while it has 
ground startup rule instances for $n \geq 2$. 

Second, note that the auxiliary and termination rules do not fully specify 
$p_{\mymax{k}}(n)$, but only its restriction to the $2\leq k \leq n$ region 
of its domain. Nevertheless, this suffices to the auxiliary purpose of the 
function in question, as it is apparent from the startup rule~(\ref{eq:5.3a}) 
and from the fact that the so restricted domain is closed under auxiliary 
reductions, these being specified by rule~(\ref{eq:5.3c}). Should a complete 
specification of $p_{\mymax{k}}(n)$ for positive integer $k$ be of interest, 
then the following two rules ought to be included as well: 
\begin{subequations}
\begin{align}
[n > 0] \qquad p_{\mymax{1}}(n) \rightarrow & \; 1 
\label{eq:5.4a}\\ 
[k > n] \qquad p_{\mymax{k}}(n) \rightarrow & \; 0 
\label{eq:5.4b} 
\end{align}
\end{subequations}

Third, as mentioned above, uniqueness of the maximal part in partitions from 
${\cal P}_{n_{\mymax{k}}}$ is warranted for $k>n/2$, thus justifying 
rule~(\ref{eq:5.3b}) for all of its domain but the boundary case $k=n/2$, 
where $n$ is even and (exactly) one of its partitions consists of two maximal 
parts, viz.\ two halves of $n$. This case ought to fall in the domain of 
rule~(\ref{eq:5.3c}), but it so happens that, precisely on these boundary 
points, the two rules turn out to be equivalent (see below). It is then 
convenient to place this part of the boundary within the domain of the 
termination rather than auxiliary rule, since this choice warrants closure 
of the aforementioned restricted domain of the subject auxiliary function 
under auxiliary construction steps, as it is argued below. 
 
Fourth, the combinatorial argument for rule~(\ref{eq:5.3c}) is similar to 
that exposed for equation~(\ref{eq:2.7}), by considering the transfer of the 
negative term to the left hand side. Then ${\cal P}_{n+1_{\mymax{k+1}}}$ may 
be split into two disjoint subsets, viz.\ the partitions where the maximal part 
has multiplicity greater than 1, and those where there's only one maximal part. 
The former are clearly counted by $p_{\mymax{k+1}}(n-k)$, by considering 
the effect of the removal of a maximal part; the latter are counted by 
$p_{\mymax{k}}(n)$, by considering the effect of subtracting 1 from the 
(only one) maximal part (which is greater than 2 by assumption, since it is 
$k+1$, with $k\geq 2$ by the domain condition; the subtraction thus does not 
make the outcome to leave the restricted domain). 

Finally, as mentioned above, rules~(\ref{eq:5.3b})~and~(\ref{eq:5.3c}) 
are equivalent for the boundary case $k\,=\,n/2$, with even $n$. Clearly, 
the previous argument applies to this special case, too, where the 
$[n\mapsto 2k]$-instance of the right hand side of rule~(\ref{eq:5.3c}) would 
be $p_{\mymax{k+1}}(2k+1) - p_{\mymax{k+1}}(k)$; the second summand here 
would vanish by rule~(\ref{eq:5.4b}), whereas by rule~(\ref{eq:5.3b}) the 
first summand would reduce to $p(k)$, which is the reduct of the 
$[n\mapsto 2k]$-instance of rule~(\ref{eq:5.3b}). 
By employing the latter rather than rule~(\ref{eq:5.3c}) for the boundary case 
in question, one thus gets the same outcome, but throws out of the latter's 
domain the only points that would lead to use of rule~(\ref{eq:5.4b}), which 
may thus be safely disposed of. 

Before embarking on the analysis of DAG construction steps produced by 
rules~(\ref{eq:5.3a}--c) in order to infer properties of the equivalent direct 
recurrence, one may note the unusual kind of induction taking place therein, 
where the auxiliary parameter increases along construction paths, until 
termination steps. This fact is easily explained by the location of the 
terminal region, which is the half region of the auxiliary recurrence domain 
that lies between the $n=k$ and the $n=2k$ boundaries (both included), 
the former coinciding with the lower boundary of the auxiliary recurrence 
domain itself. It is then fairly obvious that paths starting outside of 
the terminal region should feature increasing values of the auxiliary 
coordinate, in order to enter the terminal region eventually. 

The primary coordinate may increase as well as decrease along auxiliary paths, 
but the difference $n-k$ between the two coordinates is nonincreasing; this 
fact, together with the strictly increasing monotonicity of the auxiliary 
coordinate along construction paths warrant termination of every path 
starting at startup point $(n_0,k_0)$ after at most $n_0 - 2k_0$ 
construction steps. This can be easily seen as follows. 

The case $n_0 \leq 3 \vee k_0 \geq n/2$ is immediate, since either 
$n_0 \leq 1$, in which case there is no startup point, because a primary 
instance of rule~(\ref{eq:5.3a}) applies, or the startup point lies in the 
terminal region, hence the number of construction steps is 0. 
So, assume $n_0>3 \wedge 2\leq k_0<n/2$; then the startup point coordinates 
satisfy the domain condition of auxiliary construction rule~(\ref{eq:5.3c}), 
and its successor points along any construction path will also do so until 
fall in the terminal region. At each auxiliary construction step, along the 
path, the first coordinate distance between the source vertex of each construed 
edge and the $n=2k$ boundary of the terminal region decreases strictly, either 
by 1 if the edge is construed by the first summand in the right hand side of 
the rule (since the first coordinate increases by 2 along the aforementioned 
boundary), or by $k+2$ in the other case. The longest path up to termination 
thus consists of only edges that are construed by the first summand in the 
right hand side of the rule. Since the first coordinate distance between the 
startup point and the $n=2k$ termination boundary is $n_0 - 2k_0$, the 
similar distance between edge target vertex and the same boundary becomes 
null or negative after $n_0 - 2k_0$ construction steps at most. 

Now, about the target direct recurrence, let 
\begin{equation} 
c : \qquad n\geq 0 \quad \rightarrow \quad 
    p(n) = 1 + \mathop{\sum}_{1\leq j\leq n} c_j p(n-j) 
\label{eq:5.5}
\end{equation} 
be the equivalent direct recurrence of the composite recurrence defined by 
rules~(\ref{eq:5.3a}--c). This instance of the general form~(\ref{eq:4.1}) 
is justified by a few properties which immediately result from a first 
inspection of rules~(\ref{eq:5.3a}--c). The direct recurrence consists of 
only one recurrence equation, since there is only one startup rule and no 
primary rules in the rewriting system of the composite recurrence. The domain 
of the direct recurrence is thus as specified by the domain condition of the 
only one startup rule~(\ref{eq:5.3a}). The constant coefficient $c_{r_0}=1$ 
is also borrowed from the corresponding constant $t_r=1$ in the right hand 
side of the startup rule, since the other rules contribute null constant 
labels to the vertices of the parallel reduction DAG. 

As a matter of notation, henceforth $\tilde{n}$\ denotes the value of $n$ 
in an instance of Equation~(\ref{eq:5.5}), as well as of the primary variable 
in an instance of the startup rule~(\ref{eq:5.3a}). This is meant to prevent 
confusion with the free use of $n$ to denote the first coordinate of a 
generic point in the auxiliary plane. The DAG that is construed by 
rules~(\ref{eq:5.3a}-c) for a given $\tilde{n}$\ is then referred to as the 
$\tilde{n}$-DAG. 

Coefficients $c_j$, for $1\leq j\leq\tilde{n}$, in an instance of 
recurrence~(\ref{eq:5.5}) result from the signed unitary contributions made 
by terminating paths in the $\tilde{n}$-DAG. The sign of the contribution 
made by any given terminating path is determined by the parity of the number 
of those auxiliary construction steps in the path which lower the first 
coordinate---as the lowering corresponds to the choice of the negative literal  
in the right hand side of auxiliary rule~(\ref{eq:5.3c}), thus negative sign 
by odd parity, positive sign by even parity thereof. Regardless of sign, 
Equations~(\ref{eq:4.4}) and ~(\ref{eq:5.3b}) entail that paths in a  
$\tilde{n}$-DAG contribute to the same coefficient $c_j$ iff the coordinates 
$(n_t,k_t)$ of their terminal vertices in the auxiliary DAG 
have equal difference $n_t-k_t = \tilde{n} - j = n_0-j $, which is thus 
constant for the given $\tilde{n}$ and each given $j$, with 
$1\leq j \leq \tilde{n}$. 

\begin{figure}[ht]
\centering
\includegraphics[width=14.9cm]{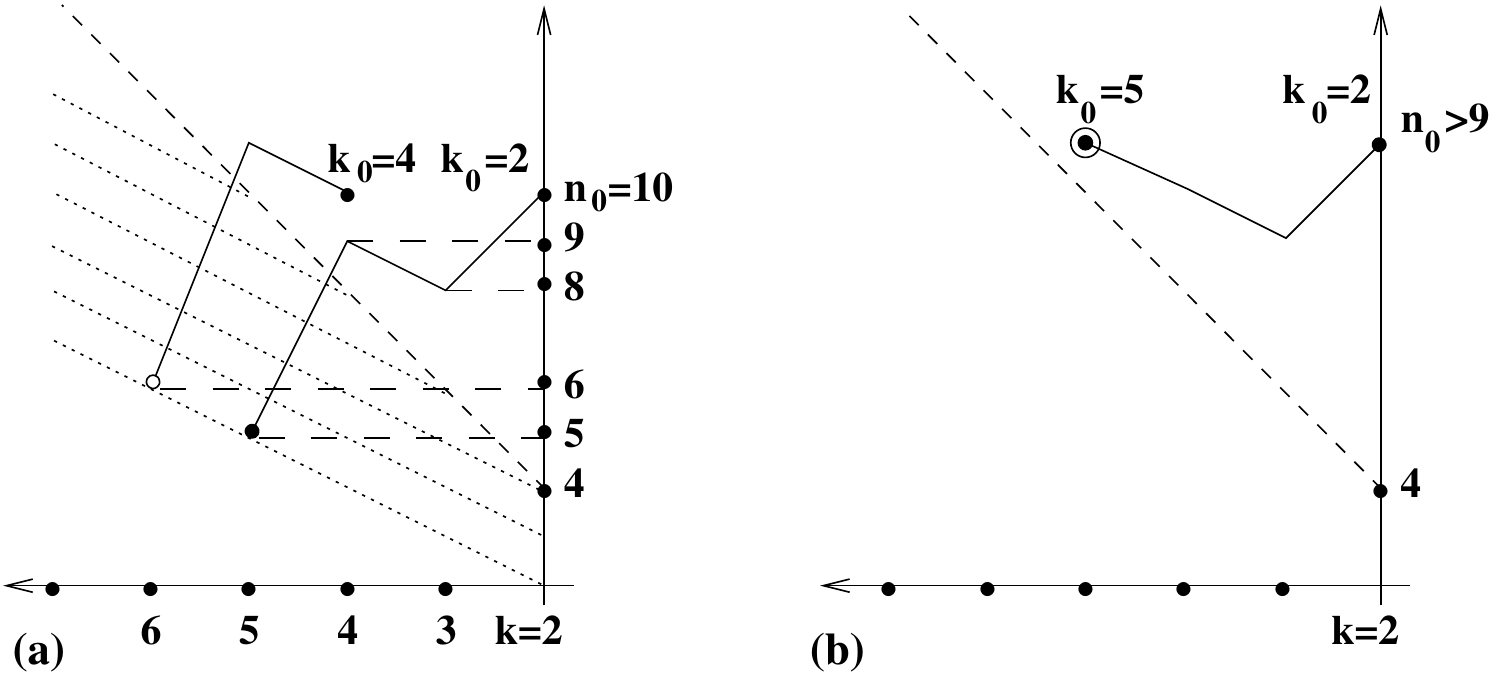}
\caption{Complementary paths by auxiliary reductions}
\label{fig:2}
\end{figure}

Recalling that the first coordinate is represented on the vertical axis, it is 
convenient in this case to give right-to-left orientation to the horizontal, 
second coordinate axis, as this choice yields a more immediate visual matching 
of the path coding which is introduced below with the graphical shape of coded 
paths. 
The domain condition of the termination rule determines the terminal region 
in the auxiliary plane, that has the diagonal line $n=k$ as lower boundary, 
the straight line $n=2k$ as upper boundary, and the $k=2$ vertical line as 
right boundary, all boundaries being included in the region. The nonvertical 
boundaries are represented in Figure~\ref{fig:2}a by the lowest dotted line 
and the broken line, where the first coordinate unit size is twice that of 
the second coordinate unit, to improve visual discrimination of paths outside 
the terminal region (thus the lowest dotted line really is the diagonal). 
The dotted lines represent equivalence classes of termination points, each 
class collecting all endpoints of terminating paths which contribute to 
the same coefficient of the target direct recurrence (not necessarily with 
the same sign). This is justified as follows. 

The aforementioned characterization of the set of paths which, regardless of 
sign, contribute to the same coefficient $c_j$ of the target recurrence, tells 
that each equivalence class is characterized by a distinct, fixed value $j$ of 
$\tilde{n}-(n_t-k_t)$, where $(n_t,k_t)$ are the termination points of the 
paths in the equivalence class. Then it is plain that, for any given 
$\tilde{n}$ and each $j$ such that $1\leq j\leq\tilde{n}$, the equivalence 
class where $n_t-k_t = \tilde{n}-j$ collects those points $(n_t,k_t)$ of the 
terminal region which lie on the straight line that is parallel to the region 
lower boundary diagonal, at vertical distance $\tilde{n}-j$ from it. 
Figure~\ref{fig:2}(a) displays two auxiliary terminating paths construed by 
the subject rewriting system for $\tilde{n}=10$, which contribute to 
$c_{10}$ with opposite signs (negative by the white-dot terminated path, 
positive by the black-dot terminated one). 

Figure~\ref{fig:2}(b) displays two paths that lie outside the terminal region 
whenever $\tilde{n}=n_0>9$, one of them consisting of a startup point 
$(n_0,5)$ only, the other starting at $(n_0,2)$ and converging to the other's 
endpoint (which happens to be the startup point in this particular example) 
after three steps. Although neither path as such contributes to any 
coefficient of the direct recurrence, it is easy to realize that every 
terminating path that starts at $(n_0,5)$, {\it i.e.}\ it consists of the 
first of the two displayed paths followed by that terminating path, may be 
put in a one-to-one correspondence with the path that starts at $(n_0,2)$ 
and proceeds as the second of the two displayed paths followed by that same 
terminating path. The two corresponding paths contribute to the same 
coefficient of the target direct recurrence, whichever that coefficient may 
be, with opposite signs, since the first displayed path has an even (null) 
number of lowering steps, whereas the second one has an odd number thereof, 
so the contributions of pairs of so correspondig paths cancel out. 

The fact that the (only) auxiliary rule~(\ref{eq:5.3c}) among the rewrite 
rules of present concern has a two-summand right hand side suggests that a 
binary representation of paths, similar to that adopted in \cite{Brown}, may 
be useful to formalizing relations and correspondences between paths. 
To this purpose, from each binary word representing a path, that will be 
referred to as the {\em path code}, it should be possible to uniquely recover 
the following information about the path: 
\begin{itemize} 
\item 
the path startup point $(n_0,x_0)$, or just $x_0$ if $n_0$ is fixed, as it is 
in our case, for each given $\tilde{n}$ determining the $\tilde{n}$-DAG of 
interest (with $n_0=\tilde{n}$ in the present case, according to the (only) 
startup rule~(\ref{eq:5.3a})); 
\item 
the sequence of binary choices between positive and negative summand in the 
auxiliary rule~(\ref{eq:5.3c}) that determines the sequence of progressively 
construed path edges.  
\end{itemize} 
The first requirement, together with the fact that, according to the startup 
rule, $k_0$ is generally unbounded entail that a subword of the path code is 
to be assigned to encoding $k_0$, whereas the rest of the code may encode the 
path edge sequence, one bit per edge. Moreover, it will be helpful to define, 
for each terminating path code, the index $j$ of the coefficient $c_j$ which 
the encoded path contributes to, and the sign of the contribution made by the 
encoded path. The former is referred to as the {\em valuation} $\nu(b)$,  
while the latter as the {\em polarity} $\pi(b)$, of argument path code $b$. 
As the example displayed in Figure~\ref{fig:2}(b) may suggest, it is actually 
useful to define these two functions on path codes for all paths in the DAG, 
regardless of whether terminating or not. 

The representation adopted in \cite{Brown} is partially fit to the present 
purpose, in the sense that, like in that case, it is convenient to let the 
indexing of bit positions in any word $b$ start from 2, this being the index 
of the rightmost bit position, thus $b=b_{l(b)+1}b_{l(b)}\ldots b_2$, with 
$l(b)$ denoting the length of word $b$. However, valuation and polarity of 
path codes differ from the corresponding definitions they are given in 
\cite{Brown}, because of the different discrete dynamics of the two composite 
recurrences under consideration. First, it is convened that the index of the 
righmost 1-bit in path code $b$ is the value of $k_0$ (consistently with the 
fact that $k_0\geq 2$ by the startup rule~(\ref{eq:5.3a})), and that each 
subsequent bit, proceeding right to left, is 1 if the corresponding edge in 
the path edge sequence is construed by the negative summand in the right hand 
side of auxiliary rule ~(\ref{eq:5.3c}), whilst it is 0 if the edge is 
construed by the positive summand in that rule. Then, letting $|b|_1$ denote 
the number of 1-bits in binary word $b$, valuation and polarity of path codes 
are defined as follows: 
\begin{subequations} 
\begin{align} 
\nu(b) & = \mathop{\sum}_{2\leq k\leq l(b)+1} kb_k , 
\label{eq:5.6a}\\
\pi(b) & = (-1)^{|b|_1+1} . 
\label{eq:5.6b} 
\end{align} 
\end{subequations} 

The polarity definition is justified by the fact that the rightmost 1-bit 
$b_{k_0}$, together with the 0-only suffix to its right, encodes $k_0$, viz.\ 
the path startup point (whose first coordinate $n_0=\tilde{n}$\ is fixed for 
all paths in the $\tilde{n}$-DAG) rather than a negative summand choice, hence 
the product of edge sign labels along the path is determined by the parity of 
the number of the {\em other} 1-bits in the path code, excluding the rightmost 
one, each of them encoding a negative sign label. 

The valuation definition meets the requirement that, if $b$ encodes a 
terminating path, then $\nu(b)$ must be the index $j$ of the coefficient 
$c_j$ which the path contributes to. For a path starting at $(n_0,k_0)$ and 
terminating at $(n_t,k_t)$, Equations~(\ref{eq:4.4}) and (\ref{eq:5.3b}) entail 
this value is $j=n_0-(n_t-k_t)=k_0+(n_0-k_0)-(n_t-k_t)$. If $t(b)$ denotes the 
number of edges in the path encoded by binary word $b$, this number clearly 
coincides with the length of $b$ minus the length of its suffix that encodes 
$k_0$, so it is given by: 
\begin{equation}
t(b) = l(b) + 1 - k_0 . 
\label{eq:5.7}
\end{equation}
Then the previous requirement on $\nu(b)$\ may be written as follows, taking 
Equation~(\ref{eq:5.7}) into account: 
\begin{equation} 
\nu(b) = k_0 + \mathop{\sum}_{1\leq i\leq t(b)} (n_{i-1}-k_{i-1})-(n_i-k_i) 
\label{eq:5.8}
\end{equation} 
The valuation $\nu(b)$\ may thus be seen to result from the sum of the initial 
second coordinate value $k_0$ and the contributions given by all edges in the 
path to bridging the gap in ``coordinate difference'' $n-k$\ between the 
startup point and the terminal point. Those edges which are construed by the 
positive summand in the right hand side of the auxiliary rule give a null 
contribution in this respect, since both coordinates increase by 1 along any 
such edge, whereas each of the other edges, say joining 
$(n_{i-1},k_{i-1})$ to $(n_i,k_i)$, contributes  
$(n_{i-1}-k_{i-1})-(n_i-k_i)=(n_{i-1}-n_i)-(k_{i-1}-k_i)%
=(n_{i-1}-(n_{i-1}-k_{i-1}))-(k_{i-1}-(k_{i-1}+1))=k_{i-1}+1=k_i$ to that 
gap reduction, according to the right hand side of the auxiliary rule. 
Now, for $1\leq i\leq t(b), k_i$\ is the index of that bit in $b$ which 
encodes the $i$-th edge in the path, therefore the contribution made by the 
$i$-th edge is $k_ib_{k_i}$\ in all cases, for $1\leq i\leq t(b)$, and since 
$k_i=k_0+i$, by Equation~(\ref{eq:5.7}) one may rewrite the contribution as 
$(k_0+i)b_{k_0+i}$ for $1\leq i\leq l(b)+1-k_0$, that is
$kb_k$ for $k_0+1\leq k\leq l(b)+1$, by an index substitution. Clearly, 
the sum of these contributions plus the initial $k_0$ is the value that 
Equation~(\ref{eq:5.6a}) assigns to $\nu(b)$, since $b_k=0$ for $2\leq k<k_0$. 
From this it immediately follows that $\nu(b)\geq 2$, hence $c_1=0$ always, 
in recurrence~(\ref{eq:5.5}). 

The double choice of the vertical axis for the first coordinate and of 
right-to-left orientation for the second coordinate axis now shows its 
comfortable effects, since it proves very easy to infer the path code from 
the visual appearance of any given path, inspected left-to-right. For example, 
the four paths displayed in Figure~\ref{fig:2} have binary codes 10100, 1011, 
1000, 0011---the reader may easily check which has which. 

A remark is in place: equation~(\ref{eq:5.6a}) defines the valuation of 
 {\em any} path code $b$, not just those of terminating paths. For a path 
with endpoint above the terminal region, this may be understood as the 
partial valuation accumulated up to that point by whichever may be the 
terminating path further proceeding from that point on. Nonetheless, it is 
useful to characterize terminating paths in terms of properties of their 
binary codes, for a given $\tilde{n}=n_0$, since only terminating paths 
deliver contributions to the target recurrence coefficients. Whether or not 
does a binary word $b$ encode a terminating path, that clearly depends 
on $\tilde{n}$, for, a change of $\tilde{n}$\ to $\tilde{n}^{_{\prime}}$ 
amounts to a vertical translation of the encoded path by 
$\tilde{n}^{_{\prime}}-\tilde{n}$, parallel to the first coordinate direction, 
and this operation on paths does not generally preserve termination. 
It will also prove useful to distinguish whether or not does the termination 
point of a terminating path belong to the upper boundary of the terminal 
region, viz.\ the straight line $n=2k$ (please note that the terminal region 
boundaries do {\em not} depend on $\tilde{n}$, as they are fixed for all DAG's; 
this fact may explain why vertical translation of paths does not preserve 
termination.) 

\begin{lemma} 
\label{lm:5.1} 
Let $b$ be the binary code of a path in the $\tilde{n}$-DAG, with 
$b_{l(b)+1}$ its leftmost bit. Then the following statements hold: 
\begin{enumerate} 
\protect\vspace*{-5mm} 
\renewcommand{\labelenumi}{(\roman{enumi})} 
\setlength{\itemsep}{1pt}
\setlength{\parskip}{0pt}
\setlength{\parsep}{0pt} 
\item 
the path encoded by $b$ is terminating iff 
$\tilde{n}-(l(b)+1)\leq\nu(b)\leq\tilde{n}$; 
\item 
the path encoded by $b$ terminates strictly below the upper boundary of the 
terminal region iff\\ 
$\tilde{n}-(l(b)+1)<\nu(b)\leq\tilde{n}$, while it terminates at that boundary 
iff $\nu(b)=\tilde{n}-(l(b)+1)$; 
\item 
if the path encoded by $b$ terminates strictly below the upper boundary of the 
terminal region, then $b_{l(b)+1}=1$; 
\item 
if $\tilde{n}-2\leq\nu(b)\leq\tilde{n}$, then $b$ encodes a terminating path 
in the $\tilde{n}$-DAG, and $b_{l(b)+1}=1$. 
\end{enumerate} 
\end{lemma} 
\noindent\textbf{Proof.}\\ 
\textit{(i)}  
The last edge of the path encoded by $b$ in the $\tilde{n}$-DAG has target 
vertex $(n_{t(b)},k_{t(b)})$ and source vertex $(n_{t(b)-1},k_{t(b)-1})$, 
where $t(b)$ is as defined by Equation~(\ref{eq:5.7}), which also entails 
$k_{t(b)}=k_0+t(b)=l(b)+1$, since $k_i$ increases by 1 at each step of the 
inductive construction of the path from its code. The path is terminating iff 
its last edge has source vertex outside the terminal region, {\it i.e.}\ 
strictly above its upper boundary, and target vertex inside it, {\it i.e.}\ 
at or below that boundary, which is the straight line $n=2k$. By the previous 
identity, this is characterized by 
$(n_{t(b)-1}>2l(b))\wedge(n_{t(b)}\leq 2(l(b)+1))$. 
The leftmost bit $b_{l(b)+1}$ makes the difference $n_{t(b)-1}-n_{t(b)}=%
b_{l(b)+1}(k_{t(b)-1}+1)-1=b_{l(b)+1}k_{t(b)}-1=b_{l(b)+1}(l(b)+1)-1$,\ 
whence the previous condition is equivalent to 
$l(b)+1\leq n_{t(b)}\leq 2(l(b)+1)$, by merging the two cases for $b_{l(b)+1}$. 
Finally, Equation~(\ref{eq:4.4}) and the stated requirement on $\nu(b)$\ give 
the identity $\nu(b)=\tilde{n}-(n_{t(b)}-k_{t(b)})$, which is equivalent to 
$n_{t(b)}=\tilde{n}-\nu(b)+l(b)+1$, whereby the previous condition proves 
equivalent to the stated one.\\ 
\textit{(ii)} 
The path encoded by $b$ in the $\tilde{n}$-DAG terminates strictly below the 
termination upper boundary iff its upward translation by 1 is a terminating 
path in the $(\tilde{n}+1)$-DAG; the replacement of $\tilde{n}$\ with 
$(\tilde{n}+1)$ in the characteristic condition provided by the previous 
statement \textit{(i)} turns the lowerbound inequality into a strict one, 
after subtraction of the added 1.\\ 
\textit{(iii)}  
If the target vertex of the last edge in the path falls strictly below the 
upper boundary of the terminal region, then that edge cannot be a rising 
one, since this would entail that also its source vertex would fall in the 
terminal region, thus outside of the auxiliary rule domain.\\ 
\textit{(iv)}  
The parallel straight lines $n-i=k$, for $0\leq i\leq 2$, that respectively 
are the path valuation equivalence classes of endpoints of those paths which 
have valuation $\tilde{n}-i$, fall entirely within the terminal region. 
Only one of these lines, viz.\ that for $i=2$, shares a point with the upper 
boundary of the terminal region, that is point (4,2), and this is the 
termination point of only one path in only one $\tilde{n}$-DAG, viz.\ the path 
encoded by $b=1$ for $\tilde{n}=4$, so $b_{l(b)+1}=1$\ holds in this case, too, 
as it does in all other subject cases by the previous statement \textit{(iii)}. 
\qed 

Finally, analogously to the binary path encoding adopted in \cite{Brown}, 
here, for every given $j\geq 2$, path codes with valuation $j$ may be put in 
a bijective correspondence with ${\cal S}_{j_{_{\!>1}}}$, the set of strict 
partitions of $j$ with smallest part greater than 1.

%% file: PentaIs6.tex
\section{Relationship with Euler's pentagonal partition and proof of the claim} 
\label{s6} 

A first bit of information about the coefficients of the direct 
recurrence~(\ref{eq:5.5}) has been easily obtained in the previous section, 
viz.\ $c_1=0$. This has a useful generalization. Recall that, according to 
\cite{Z2008}, a function $a_n$ is said {\it C-recursive} if it satisfies a 
linear recurrence with constant coefficients, 
$a_n = c_1 a_{n-1}+\ldots+c_d a_{n-d}$. 

\begin{lemma} 
\label{lm:6.1} 
The recurrence~(\ref{eq:5.5}) has constant coefficients, {\em i.e.}\ $c_j$ 
only depends on $j$, not on $n$. 
\end{lemma} 
\noindent\textbf{Proof.} 
For every $\tilde{n}>2$, a bijection is established between terminating paths 
that have the same polarity and the same valuation $j<\tilde{n}$\ in the 
$\tilde{n}$-DAG and in the $(\tilde{n}-1)$-DAG. 
Every path in the $(\tilde{n}-1)$-DAG that terminates strictly below the upper 
boundary of the terminal region, viz.\ the $n=2k$\ straight line, is mapped 
to the path in the $\tilde{n}$-DAG that has the same binary code; this map is 
clearly injective, and it amounts to an upward-by-1 translation of paths from 
the $(\tilde{n}-1)$-DAG into the $\tilde{n}$-DAG, along the first coordinate 
direction. 
The same mapping rule would not work for paths in the $(\tilde{n}-1)$-DAG that 
terminate at the upper boundary of the terminal region, since the image 
path under translation would not be a terminating path in the $\tilde{n}$-DAG. 
Therefore, for every path in the $(\tilde{n}-1)$-DAG that terminates at 
the upper boundary and has binary code $b$, its bijective image in the 
$\tilde{n}$-DAG is the path which has binary code $0b$, that is easily seen 
to be a terminating one, also at the upper boundary of the terminal region. 
The so defined map is a bijection, thanks to Lemma~(\ref{lm:5.1}(iii)), which 
entails disjointness of the images of the aforementioned two classes of 
terminating paths in the $(\tilde{n}-1)$-DAG under the respective mapping rules 
as given above. This bijection includes all terminating paths with valuation 
$j<\tilde{n}$, and it preserves both valuation and polarity, so the resulting 
value of $c_j$ is the same in both DAG's, for all $j<\tilde{n}$. 
\qed 

Thanks to Lemma~(\ref{lm:6.1}), it suffices to compute each $c_j$ for the 
smallest $\tilde{n}$ for which $c_j$ is defined, that is $\tilde{n}=j$, 
since it will thereafter keep constant for all higher values of $\tilde{n}$. 
This fact leads to an almost surprisingly simple proof of 
Equation~(\ref{eq:1.4}). 
Two more lemmas provide useful tools to that purpose. 

\begin{lemma} 
\label{lm:6.2} 
The following identity holds for all $j\geq 1$: 
\[c_j = \mathop{\sum}_{\nu(1b)=j} \pi(1b)\] 
\end{lemma} 
\noindent\textbf{Proof.} 
The sum is null for $j=1$, consistently with the already assessed $c_1=0$. 
For $j>1$, by Lemma~(\ref{lm:6.1}) it suffices to compute $c_j$ in the 
$\tilde{n}$-DAG where $\tilde{n}=j$. Lemma~(\ref{lm:5.1}(iv)) tells that 
all paths which have valuation $j$ are terminating paths in the $j$-DAG and 
have leftmost bit 1 in their path code. 
\qed 

The paths in the $j$-DAG that have valuation $j$ are those which terminate at 
the lower boundary of the terminal region, viz.\ the diagonal line $n=k$. 
The next statement is a useful tool to compute the valuation of a path code 
out of the valuation of subwords of its. 
\begin{lemma} 
\label{lm:6.3} 
If $b_p,b_s$ are binary words and their concatenation is $b = b_p b_s$, then 
\begin{equation} 
\nonumber 
\nu(b) = \nu(b_s) + \nu(b_p) + |b_p|_1 l(b_s) 
\end{equation} 
\end{lemma} 
\noindent\textbf{Proof.} 
Follows from the definition~(\ref{eq:5.6a}) of the valuation function. 
\qed 

The previous lemmas are all that is needed to show validity of the main claim.

\subsection{Bijective proof} 
\label{s6.1} 

\begin{proposition} 
The coefficients $c_j$ of recurrence~(\ref{eq:5.5}) satisfy $c_j=f_j$ for all 
$j>0$, with $f_j$ defined by Equation~(\ref{eq:1.3}). 
\label{pr:6.1} 
\end{proposition} 
\noindent\textbf{Proof.} 
By induction on $j$. The basis case $j=1$ is immediate, since $e_0+e_1=0$ by 
Equation~(\ref{eq:1.1}). For the inductive step, it suffices to show that 
$c_j-c_{j-1}=f_j-f_{j-1}$, thanks to the induction hypothesis. Since 
$f_j-f_{j-1}=e_j$ by Equation~(\ref{eq:1.3}), then by Lemma~(\ref{lm:6.2}) 
it suffices to find an involution $\leftrightarrow$\ on the set of binary path 
codes $B_j\cup B_{j-1}$, with $B_i\meqdef\{1b \; | \; \nu(1b)=i\}$, that 
satisfies the following requirements: 
(i) if $1b\leftrightarrow 1b^{\prime}$ and $b\neq b^{\prime}$, 
then $\pi(1b)=\pi(1b^{\prime})$\ iff $\nu(1b)\neq\nu(1b^{\prime})$, and 
(ii) $1b\leftrightarrow 1b$\ iff $\nu(1b)=j$\ is pentagonal, 
say $j=(3k^2\pm k)/2$, in which case $\pi(1b)=(-1)^{k+1}$. 
Such an involution may be specified by as few as two {\em mapping rules}, 
which are pairs of {\em word patterns}; these are words over the binary 
alphabet extended with variables which range over binary words such that 
the pattern instance meets a specified {\em domain condition}. 
The first rule to this purpose defines a bijection between $B_{j-1}$ and the 
subset of $B_j$ that consists of those path codes which have a 10 prefix; 
putting brackets around domain conditions, and using the abbreviation 
``$[D]\;w\leftrightarrow w^{\prime}$'' to stand for 
``$[w\!\in\!D]\;w\leftrightarrow w^{\prime}$'', here is this rule: 
$[B_j]\;10x \leftrightarrow 1x$. Please note that the specified domain 
condition, which applies to the binary word instances of the left hand side 
pattern, together with Lemma~\ref{lm:6.3} entail the converse domain condition 
$[1x\!\in\!B_{j-1}]$ for the binary word instances of the right hand side 
pattern. The second mapping rule is actually a rule scheme, since the set of 
its constituents is extended with a variable ranging over the nonnegative 
integers, subject to validity of the domain condition. Here it is: 
$[B_j]\;1^{k+2}0x0^{k+1} \leftrightarrow 1^{k+2}x10^k$. 
Finally, the fixed points of the involution are defined as those words in 
$B_j$ to which neither rule assigns a correspondent. The rest of the proof 
consists of a straightforward check of the following facts. 
(1) The converse domain condition for the second rule is 
$[1^{k+2}x10^k\!\in\!B_j]$ (thus corresponding binary instances of the two 
word patterns get the same valuation, viz.\ $j$). 
(2) Corresponding binary instances of the first rule have the same polarity, 
whereas those of the second rule have opposite polarity. 
(3) The previous two facts entail validity of requirement (i). 
(4) Fixed points of the involution are all the binary words in $B_j$ that take 
any of the following forms, for $k\geq 0$: $1^{k+2}0^k, 1^{k+2}0^{k+1}, 1$. 
(5) Valuations of these fixed point path codes respectively are the pentagonals 
$(k+2)(3(k+2)-1)/2, (k+2)(3(k+2)+1)/2, 2$; the only pentagonal that is not 
captured by any of these forms is 1, but this falls outside of the $B_j$\ part 
of any domain $B_j\cup B_{j-1}$ of the subject involutions, since 2 is the 
smallest valuation $j$ of present concern, nor is it relevant to the inductive 
step of the proof, since $j=1$ is the basis case; the first clause of 
requirement (ii) is thus satisfied. 
(6) Polarities of the two families of fixed point path codes are 
$\pi(1^{k+2}0^k)=\pi(1^{k+2}0^{k+1})=(-1)^{k+3}$ by Equation~(\ref{eq:5.6b}), 
hence $(-1)^{(k+2)+1}$, thus fulfilling the last clause of requirement (ii) 
for $j>2$. 
(7) Fixed point $1$ also meets the last clause of requirement (ii) for the 
$j=2$ case, with positive polarity by Equation~(\ref{eq:5.6b}), which fact 
completes the proof. 
\qed 

A final remark about language-theoretic sideways of the previous proof may 
be of interest to some readers. It seems that a key factor behind the great 
parsimony in the number of mapping rules that suffice to formalize the 
involution in the previous proof, is the particular selection of binary codes 
which have pentagonal valuations and that form the set of fixed points. 
A quick look at the pattern of the two families, with the nonnegative integer 
$k$ as pattern variable, tells that they do not form a regular language over 
the binary alphabet, rather a context-free one, whose path codes may be 
visualized as the ``trapezoidal'' Ferrers diagrams, displayed {\it e.g.}\ in 
\cite{Z2003}, which play a key r\^ole in Franklin's proof \cite{F1881}. This 
fact is easily realized by taking the remark at the end of Section~\ref{s5} 
into account. 
One may well take a different set of binary codes as representatives of the 
pentagonal numbers, that does form a regular language. The following regular 
expression testifies to this possibility: 
$(100)^*(1\,${\footnotesize$\,+\,$}$\,011)$ (with ``$\,^*\,$'' and 
``{\footnotesize $\,+\,$}'' the regular Kleene star and choice operators, 
respectively). However, the author must admit his proven inability to build 
an involution that would isolate these path codes as fixed points.

\subsection{Proof by generating functions} 
\label{s6.2} 
As pointed out at the end of Section~\ref{s5}, terminating path codes with 
valuation $n$ may be put in a one-to-one correspondence with the strict 
partitions of $n$ that have smallest part greater than 1. Essentially, this 
means that the indices of 1-bits in the path code are the (necessarily 
distinct) parts in the corresponding strict partition of the valuation of the 
path code itself. The polarity of the path code thus uniquely corresponds to 
the parity of the number of parts in the corresponding partition, odd parity 
corresponding to positive polarity. According to Lemma~(\ref{lm:6.2}), 
coefficient $c_j$ thus results from the difference $O(j)-E(j)$ between the 
number of strict partitions of $j$ with an odd number of parts and that of 
such partitions with an even number of parts, all partitions being constrained 
to have smallest part greater than 1. It takes a little effort of combinatorial 
imagination to identify the following generating function as that which suits 
the present purpose: 
\begin{equation} 
\label{eq:6.1} 
-\mathop{\prod}_{j\geq 2}(1-x^j) = \mathop{\sum}_{n\geq 0}c_n\,x^n , 
\end{equation} 
the negative sign on the left hand side being explained by the fact that 
selecting the $x^j$ term in an odd number of binomials must yield a positive 
contribution to the relevant coefficient on the right hand side, and 
conversely for selection of an even number of $x^j$ terms. 
From Equation~(\ref{eq:6.1}) we may immediately infer $c_0=-1$ and $c_1=0$. 
The following manipulation of Equation~(\ref{eq:6.1}) showcases a general 
method of getting recurrences out of generating functions \cite{W2000}. 

Let $F\meqdef\mathop{\sum}_{n\geq 0}c_n\,x^n$. Taking logarithms in 
Equation~(\ref{eq:6.1}), then turning the left hand side into a sum, and 
finally taking derivatives yields the following identity: 
\begin{equation} 
\label{eq:6.1a} 
\nonumber 
\mathop{\sum}_{j\geq 2}\frac{-jx^{j-1}}{1-x^j}F = 
F^{\prime} = 
\mathop{\sum}_{n\geq 0}nc_nx^{n-1} . 
\end{equation} 
The following identity is then worked out, where $\frac{1}{1-x^j}$ is 
replaced with the geometric series $\mathop{\sum}_{k\geq 0}x^{jk}$: 
\begin{equation} 
\label{eq:6.1b} 
\nonumber 
\frac{-jx^{j-1}}{1-x^j} = \frac{j}{x}\frac{-x^j}{1-x^j} = 
\frac{j}{x}(1-\frac{1}{1-x^j}) = -\frac{j}{x}\mathop{\sum}_{k\geq 1}x^{jk} = 
-j\mathop{\sum}_{k\geq 1}x^{jk-1} . 
\end{equation} 
By introducing the right hand side of this equation into the previous one, and 
therein expanding $F$, with a renaming of its index for clarity of later 
manipulation, one gets the following: 
\begin{equation} 
\label{eq:6.1c} 
\nonumber 
-\left(\mathop{\sum}_{j\geq 2}j\mathop{\sum}_{k\geq 1}x^{jk-1}\right)%
\left(\mathop{\sum}_{i\geq 0}c_i\,x^i\right) 
= 
\mathop{\sum}_{n\geq 0}nc_nx^{n-1} . 
\end{equation} 
By equating the coefficients of $x^{n-1}$ on both sides one then gets: 
\begin{equation}
\label{eq:6.1d} 
\nonumber 
nc_n = -\mathop{\sum}_{0\leq i\leq n\bminus 2}c_i\mathop{\sum}_{j\geq 2}j%
\mathop{\sum}_{\substack{k\geq 1\\ jk\bminus 1=n\bminus 1\bminus i}}1 . 
\end{equation} 
Now, the condition $jk-1=n-1-i$ is satisfied iff $j|(n-i)$, and for each such 
$j$ there is a unique $k=\frac{n-i}{j}$\ fit to the purpose, therefore the 
inner double summation may be equated to $\sigma(n-1)-1$, where $\sigma$ is 
the sum of divisors function, and the outer -1 is due to the exclusion of 
$j=1$ from the count of the divisors of $n-1$, since $j\geq 2$ is required 
by the second summation indexing. 
One finally gets the following recurrence for the coefficients $c_n$ specified 
by the generating function ~(\ref{eq:6.1}), {\it i.e.}\ the coefficients of 
the target recurrence~(\ref{eq:5.5}): 
\begin{equation} 
\label{eq:6.2} 
c_n = -\frac{1}{n}\mathop{\sum}_{0\leq i\leq n\bminus 2}(\sigma(n-i)-1)c_i 
\end{equation} 
The similar manipulation of the well-known generating function~\cite{E1783} 
for the coefficients of Euler's pentagonal recurrence~(\ref{eq:1.2}) 
yields the following recurrence for them, with basis $e_0=-1$: 
\begin{equation} 
\label{eq:6.3} 
e_n = -\frac{1}{n}\mathop{\sum}_{0\leq i\leq n\bminus 1}\sigma(n-i)e_i 
\end{equation} 
By Equation~~\ref{eq:1.3}, the following proposition is clearly equivalent 
to Proposition~\ref{pr:6.1}, but the proof exploits the recurrences obtained 
from the respective generating functions for the subject coefficients. 
\begin{proposition} 
The coefficients $c_j$ of recurrence~(\ref{eq:5.5}) satisfy $c_j-c_{j-1}=e_j$,\ 
for all $j>0$. 
\label{pr:6.2} 
\end{proposition} 
\noindent\textbf{Proof.} 
Equation~(\ref{eq:6.1}) gives $c_0=-1$ and $c_1=0$; these identities show 
validity of the basis in the proof of the statement by induction on $j$, viz.\ 
for $j=1$; the inductive step follows by manipulating Equations~(\ref{eq:6.2}) 
and (\ref{eq:6.3}), using the induction hypothesis (IH). Here are the main 
steps, with concise justifications in brackets, the reader should be able to 
fill the gaps. Assume $c_i-c_{i-1}=e_i$\ for $0<i\leq j$\ as IH, then rewrite 
$e_{j+1}$\ as follows:   
\begin{align*} 
[Eq.~(\ref{eq:6.3})] \; 
e_{j+1} &= -\frac{1}{j+1}\mathop{\sum}_{0\leq i\leq j}\sigma(j+1-i)e_i\\ 
[IH] \qquad \qquad \;\, 
&= -\frac{1}{j+1}\left(\sigma(j+1)c_0%
   +\mathop{\sum}_{1\leq i\leq j}\sigma(j+1-i)(c_i-c_{i-1})\right)\\ 
[Eq.~(\ref{eq:6.2})] \qquad 
&= -\frac{1}{j+1}\left(-(j+1)c_{j+1}+\left(\mathop{\sum}%
_{0\leq i\leq j\bminus 1}c_i\right)+\sigma(1)(c_j-c_{j-1})\right)\\ 
& \quad -\frac{1}{j+1}\left(\left(-\mathop{\sum}_{1\leq i\leq j\bminus 1}%
    (\sigma(j-(i-1))-1)c_{i-1}\right)%
    +\mathop{\sum}_{1\leq i\leq j\bminus 1}c_{i\bminus 1}\right)\\ 
[Eq.~(\ref{eq:6.2})] \qquad 
&= -\frac{1}{j+1}\left(-(j+1)c_{j+1}%
   +\left(\mathop{\sum}_{0\leq i\leq j-1}c_i\right)%
   +\sigma(1)(c_j-c_{j-1})+jc_j-\mathop{\sum}_{0\leq i\leq j-2}c_i\right)\\ 
&= c_{j+1}-\frac{j}{j+1}c_j-\frac{1}{j+1}(c_{j\bminus 1}+\sigma(1)%
(c_j-c_{j\bminus1}))\\ 
&= c_{j+1}-c_j 
\end{align*} 
\qed

%% file: PentaIs7.tex
\section{Conclusions} 
\label{s7} 
Neither novelty nor computational efficiency justify interest in the 
recurrence for integer partition investigated in this work. On the novelty 
side, as pointed out to the author by Nick Loehr\cite{L2010}, the subject 
recurrence is essentially that proposed in Exercise 5.2.3 of Igor Pak's 
survey \cite{P2006}, although it is not noted there that the coefficients 
result from a discrete integration of Euler's coefficients; so, the essence 
is the same, the form is different, and a new form may raise some interest 
at times. On the computational side, the present recurrence, albeit linear 
and C-recursive, is less efficient than Euler's pentagonal recurrence, as 
the latter requires the computation of fewer recurrents. 
What seems more interesting is the sort of duality between the respective 
(bijective) proof techniques which extract them from composite recurrences. 
Euler's recurrence may be termed the ``derivative'' pentagonal recurrence, 
and may be obtained by induction on minimal parts; the recurrence presented 
here might be termed the ``integral'' pentagonal recurrence, and is obtained 
by induction on maximal parts. Is this a situation which is peculiar to the 
integer partition function, or does it occur in other situations? Should the 
latter be the case, under which general conditions ought it to be expected? 

Another aspect which may be of some interest is the fact that, while less 
efficient on the computational side, the ``integral'' pentagonal recurrence 
is obtained by what seems to be a more parsimonious construction of the 
bijection, if one compares the bijection rules presented here with those 
worked out in \cite{Brown}, which is the closest case to carry out such a 
comparison. A possibly interesting aside of this observation is that, if one 
replaces Equation~(\ref{eq:1.3}) with its ``derivative'' counterpart, viz.\ 
$e_0=-1, e_i=f_i-f_{i-1}$, taken as a definition of Euler's coefficients, then 
the bijective proof here presented for the ``integral'' pentagonal recurrence, 
 {\em together with} (an aptly rearranged variant of) the proof by generating 
functions yield a novel proof of the well-known fact that Euler's coefficients 
are recurrence coefficients for integer partition. We invite the reader to try 
to show this fact as an exercise. 

Another intriguing question, accompanied by a probably more challenging kind 
of exercise for the curious reader, is posed at the end of Section \ref{s6.1}. 
Such kind of questions naturally arise in the context of bijective proofs, 
with some finitary language encoding of the set whereon an involution is 
sought for. Language theoretic questions and approaches enjoy some popularity 
in algebraic combinatorics, see {\em e.g.}\ the recent \cite{CRZ2010} for an 
exciting new perspective on a century-old problem. 

The most relevant contribution made by this note to the author's own research 
interest is in the {\em method} adopted, to turn composite recurrences into 
direct ones, under relatively mild assumptions. It is {\em not} an automated 
method, but it seems to support combinatorial reasoning and to prove helpful 
to combine arguments of different kinds, e.g. bijective vs.\ generating 
functions, in the construction of proofs of equivalence of different 
recurrences as well as in the discovery of new, direct recurrences. 
As a matter of fact, that's how the present result, which the title of this 
note is about, came to the fore; prompted by the intriguing statement in 
\cite{Brown}, that taking the number of partitions of $n$ with the smallest 
term $j$ is {\em one} way of approaching the problem, the problem being to get 
a better understanding of why Euler's Pentagonal Theorem is true, it seemed 
just natural to try the dual way, viz.\ that of taking the number of partitions 
of $n$ with the largest term $j$. The aim was to find another proof of Euler's 
recurrence, but the surprising outcome was a different, equivalent recurrence. 

Notwithstanding the author's excitement about the method showcased in this 
note, its actual value is far from being assessed. While it seems reasonable 
to expect to find it useful with linear, C-recursive recurrences, more complex 
recurrence kinds may give rise to new challenges. This will be a subject of 
further investigation in the near future.

%% file: PentaIak.tex
\section*{Acknowledgements} 

The author wishes to thank Vincenzo Manca for his gentle introduction to 
the exciting world of Algebraic Combinatorics and Partition Theory, and for 
the timely provision of those references which gave birth to the compelling 
need to write the present note.